\documentclass[12pt]{amsart}

\title{More about additive representation functions for integers} 

\author{Labib HADDAD}
\address{120 rue de Charonne, 75011 Paris, France}
\email{labib.haddad@wanadoo.fr}

\usepackage{amssymb}
\usepackage{amsmath}
\usepackage{endnotes}

\newcommand{\su}{\subsection*}

\newcommand{\head}{\section*}
\newcommand{\noi}{\noindent}

\newcommand{\leqs}{\leqslant}

\newcommand{\geqs}{\geqslant}

\newcommand{\N}{\mathbb N}

\begin{document}
\maketitle

\thispagestyle{empty}

\markboth{Labib Haddad}{Representation functions}

\

In a recent  note [3], posted on {\tt arXiv}, (16 Jul 2015),   {\sc Kiss} and {\sc 
S\'andor}  \lq\lq \sl improve a result of Haddad and Helou  about the Erd\"os-Tur\'an conjecture"\rm. Can this improvement be still improved ? We try to answer that question. 

\ 

Here are a few definitions and notations, stripped down as much as can be,   for simplicity sake.

\

Let $X$ be any subset of $\N = \{0, 1, \dots\}$. A $X$-{\bf representation} of an integer $n$ is any ordered couple $(x,y)\in X\times X$ such that $n = x+y$. Let $r(X,n)$ be the number of all those $X$-representations of $n$. The function  $n\mapsto r(X,n)$ is the {\bf representation function} relative to $X$.  Set $s(X) = \sup \{r(X,n) : n\in \N\}$.  Now, $s(X)$ is either an integer or $\infty$, according to cases, and it produces, thus, a {\bf dichotomy} among the subsets of $\N$. We say that $X$ is in {\bf the upper class} whenever $s(X) = \infty$. \ Otherwise, say that $X$ is in {\bf the lower class}. 

\

\noi All along, $A$ and $B$ are two given infinite subsets of $\N$.  We enumerate $A$ and $B$ in increasing order:
$$A = \{a_1 < a_2< \dots\} \ \ , \ \ B = \{b_1 < b_2 < \dots\},$$
then set
$$A(k) = \{a_1 < a_2< \dots < a_k\} \ \ , \ \ B(k) = \{b_1 < b_2 < \dots < b_k\},$$
$$u(k) = s(A(k)) \ , \ v(k) = s(B(k)).$$
Clearly enough, $u(k)$ and $v(k)$ are monotone, non-decreasing, functions of $k$, and $s(A)$, $s(B)$, are their respective limits when $k$ approaches $\infty$.
Also, set:
$$d(k) = \sup \ \{|a_i-b_i| : i\leqs k\},$$
$$d = \sup\{d(k) : k\geqs 1\} = \lim_{k\to\infty} d(k).$$
Thus, $d$ is a measure of proximity, or closeness, between the two subsets $A$ and $B$. It seems reasonable to expect that whenever $A$ and $B$ are close enough, in a certain sense, they both belong to the same class, upper or lower. That this is indeed the case has been already noticed a long time ago. More specifically,  {\bf if $\mathbf d$ is finite, then $\mathbf A$ and $\mathbf B$ belong, both,  to the same class}: See, for example, in [1], Corollary 3.4, page 88, and the inequalities:
\[\frac{s(A)}{4d + 1}  \leqs s(B) \leqs (4d + 1)s(A).\tag{0}\]

\

{\bf What if $\mathbf d$ is infinite ?}

\

\noi Of course, for $k\geqs 1$, the following more general inequalities  still hold:
\[\frac{u(k)}{4d(k) + 1} \leqs v(k) \leqs (4d(k) + 1)u(k).\tag{1}\]
\su{Proof} For $m, n\in\N$, set
$$E(k,m) = \{(i,j) : i, j \leqs k \ , \ a_i + a_j = m \},$$
$$F(k,n) = \{(i,j) : i, j \leqs k \ , \ b_i + b_j = n \}.$$
For $(i,j) \in F(k,n)$, we have
$$b_i - d(k) \leqs a_i \leqs b_i + d(k)$$
$$b_j - d(k) \leqs a_j \leqs b_j + d(k)$$
so that
$$b_i + b_j - 2d(k) \leqs a_i + a_j \leqs b_i + b_j + 2d(k),$$
$$n - 2d(k) \leqs a_i + a_j \leqs n + 2d(x).$$
Therefore, each couple $(i,j)\in F(k,n)$  belongs to one of the sets $E(k,m)$ for some $m \in [n-2d(k),n+2d(k)]$.
The number of couples in $F(k,n)$ is $r(B(k),n)$. The number of couples in $E(k,m)$ is $r(A(k),m)\leqs u(k)$. Taking for $n$ an integer having the maximum number of $B(k)$-representations, i.e., $r(B(k),n) = v(k)$,  we thus obtain:
 $$v(k) \leqs (4d(k) + 1)u(k).$$
Exchanging $A$ and $B$ obtains the result.\qed

\

\noi Let us introduce a new function:
\[w(A,B) = \sup_{k\geqs 1}\frac{u(k)}{4d(k)+1}.\tag{2}\]
\{One might as well call $w$ the {\it wizard}.\}
Then, clearly enough, we have 
\[w(A,B)\leqs s(B) \ \text{and} \ w(A,B) \leqs \sup u(k) = s(A).\tag{3}\]  We similarly have, of course, $w(B,A) \leqs s(A)\land s(B)$. 
\su{Scholium} \sl Two given subsets $A$ and $B$ of $\N$ are both in the upper class  if $w(A,B)$ is infinite\rm. Indeed, if $w(A,B)$ is infinite, then so are $s(A)$ and $s(B)$, by (3).

\

\head{A special case}

\

Take $A$ to be the set of squares, $A = \{1,4,9,\dots, n^2, \dots\}$. Consider $u(k) = s(A(k))$. It is well-known that $u(k)$ is unbounded, that is, $s(A)$ is infinite: Just remember, for instance, Jacobi's formula for the number of representations of an integer as a sum of two squares. 

\

\noi Take any function $f(k)>0$  such that $u(k)/f(k)$ is unbounded. If, for a given subset $B$ of $\N$, we have $d(k)\leqs f(k)$, then $B$ is in the upper class.

\

\noi Otherwise stated: {\bf Whenever we have $\mathbf{|b_n - n^2| \leqs f(k)}$, for each $n\leqs k$, the subset $\mathbf B$ belongs to the upper class.}

\

\noi One way to choose such a function $f$ is to take a function $g(k) >0$ such that $\lim_{k\to\infty}g(k) = 0$, and let $f(k) = u(k)g(k)$. Examples abound.

\

\noi This is an ample generalization of theorem 2 in [3]. 

\

\noi \{See [2] for the mentioned result of Haddad and Helou.\}

\

\head{References}

\

1.  G. Grekos, L. Haddad, C. Helou, and J.  Pihko, \sl The class of Erd\"os-Tur\'an sets\rm,  Acta Arith. {\bf 117} (2005), 81-105. 

\

2.  L. Haddad, C. Helou,  \sl Representations of integers by near quadratic sequences\rm , Journal of Integer sequences {\bf 15} (2012), 12.8.8.

\

3.  S\'andor Z. Kiss and Csaba S\'andor, \sl On the maximum values of the additive representation functions", {\tt arXiv:1504.07411v2}. 

\

\

\enddocument